\documentclass[11pt]{amsart}
\usepackage{txfonts}
\usepackage{amssymb}
\theoremstyle{plain}

\newtheorem{Thm}{Theorem}

\newtheorem{Def}[Thm]{Definition}

\begin{document}

\title[Yamabe metrics, Fine solutions to the Yamabe flow ]
{Yamabe metrics, Fine solutions to the Yamabe flow, and local $L^1$-stability}

\author{Li Ma}
\address{Li MA,  School of Mathematics and Physics\\
  University of Science and Technology Beijing \\
  30 Xueyuan Road, Haidian District
  Beijing, 100083\\
  P.R. China }
\address{  Department of Mathematics \\
Henan Normal university, Xinxiang, 453007 \\
China}
 \email{lma@tsinghua.edu.cn}

\thanks{Li Ma's research is partially supported by the National Natural
  Science Foundation of China (No.11771124)}

\begin{abstract}
In this paper, we study the existence of complete Yamabe metric with zero scalar curvature on an n-dimensional complete Riemannian manifold $(M,g_0)$, $n\geq 3$. Under suitable conditions about the initial metric, we show that there is a global fine solution to the Yamabe flow. The interesting point here is that we have no curvature assumption about the initial metric.
We show that on an n-dimensional complete Riemannian manifold $(M,g_0)$ with non-negative Ricci curvature, $n\geq 3$,
the Yamabe flow enjoys the local $L^1$-stability property from the view-point of the porous media equation.
Complete Yamabe metrics with zero scalar curvature on an n-dimensional Riemannian model space are also discussed.


\end{abstract}

\maketitle

\section{Introduction}\label{sect0}

In this paper, we study several problems about the Yamabe metrics with zero scalar curvature on complete non-compact Riemannian manifolds. One problem under consideration is to study the existence of the Yamabe metric with zero scalar curvature on the complete Riemannian manifold $(M,g_0)$ with non-negative scalar curvature. We develop a neat linear theory for this problem.
We consider the fine solutions (see below for the definition) of Yamabe flow near the known Yamabe metric with zero scalar curvature on the complete Riemannian manifold $(M,g_0)$ with non-negative Ricci curvature. We also consider the local $L^1$ stability and uniqueness of the Yamabe flow on such a manifold from the view-point of the porous media equation \cite{HP}. We also discuss the Yamabe metric with zero scalar curvature on radially symmetric manifolds.

We now recall some facts about the Yamabe metric on the complete manifold $(M^n,g_0)$, $n\geq 3$.
We define the conformal Laplacian operator by
$$
Lu=-\frac{4(n-1)}{n-2} \Delta_{g_0} u+R_0 u, \ \ 0<u\in C^2(M).
$$
 Let $g=u^{\frac{4}{n-2}}g_0$. Then the scalar curvature of the metric $g$ is
 $$
 R(g)=u^{-\frac{n+2}{n-2}}Lu.
 $$
The largely open question is to find a complete Yamabe metric on $(M^n,g_0)$, that is, the complete metric $g$, conformal to $g_0$, with the scalar curvature $R(g)=$constant.  This is called the Yamabe problem on complete non-compact Riemannian manifolds (see problem 32 in page 677 in the book \cite{SY} and there is a counter-example given by Z.Jin \cite{J} in the sense that there is a lot complete non-compact Riemannian manifolds which have no complete Yamabe metric conformal to the given ones.

 The difficulty in understanding the Yamabe problem on complete non-compact Riemannian manifolds is the lack of compactness and there are not much methods to handle it.
 If $(M,g_0)$ is compact, R.Hamilton \cite{H} has proposed the Yamabe flow to get the Yamabe metric. This flow method is successful and there are interesting progresses made by B.Chow \cite{C}, R.Ye \cite{Y}, H.Schwetlick and M.Struwe \cite{SS}, and S.Brendle \cite{B}. When $(M,g_0)$ is a complete non-compact Riemannian manifold, we may set up the Yamabe flow by the domain exhaustion method \cite{MA}.
 Recall that the Yamabe flow is the metric family $(g(t))$ ($t\in [0,T$, $T>0$)  such that
\begin{equation}\label{flow}
\partial_t g= -R(g)g, \ \  g=g(t)=u(t)^{\frac{4}{n-2}}g_0, \ \ t\in [0,T],
\end{equation}
with the given initial metric $g(0)$. AS one may expect, when the initial metric has non-negative scalar curvature, the limit of the Yamabe flow may tend to a complete Riemannian metric with zero scalar curvature \cite{M} and \cite{M1}. Yamabe flow is an important tool to understanding the geometric structure for pinched metrics \cite{MC}, which a Myers type theorem has been obtained by the method of Yamabe flow. We say a global solution to \eqref{flow}  a \emph{fine solution} to \eqref{flow} if there are uniform positive constants $C_1$ and $C_2$ such that
$$
C_1g_0\leq g(t)\leq C_2g_0, \ \ t\in [0,T].
$$
That is to say, the metric $g(t)$ is uniformly quasi-isometric equivalent to the given metric $g_0$.
For $p\in M$, we denote by $d_{g_0}(x,p)$ the distance between the two point $x\in M$ and $p\in M$ in the metric $g_0$. Sometimes we may write $|x|=d(x,p):=d_{g_0}(x,p)$ for simplicity.
We shall prove below that, if $(M,g_0)$ is a complete Riemannian manifold with non-negative scalar curvature and with Ricci curvature bounded from below and if the Poisson equation $$
-\frac{4(n-1)}{n-2} \Delta_{g_0} u+R_0u=R_0, \ \ in \ \ M, \ \ u(x)\to 0\ \ as \ d(x,p)\to \infty,
$$
has a positive solution $v$, then for a large class of initial metrics, there are global fine solutions to the Yamabe flow.
We remark that we have no curvature assumption about the initial metric $g(0)$ except $g(0)$ is quasi-isometric equivalent to the given metric $g_0$. Roughly speaking, the idea is below.
Let $C\geq 1$ be a positive constant and consider $v_-=1-Cv$. Then
$$
-\frac{4(n-1)}{n-2} \Delta_{g_0} v_-+R_0v_+=R_0(1-C)\leq 0, \ \ on \ \ M,
$$
which means that $v_-$ is a sub-solution to the equation
\begin{equation}\label{yamabe0}
-\frac{4(n-1)}{n-2} \Delta_{g_0} u+R_0u=0, \ \ on \ \ M.
\end{equation}
Similarly, letting $v_+=1+Cv$, we have
$$
-\frac{4(n-1)}{n-2} \Delta_{g_0} v_++R_0v_+=R_0(1+C)\geq 0, \ \ on \ \ M.
$$
That is to say, $v_+$ is a super solution to the equation \eqref{yamabe0}. Then, by the monotone method we have a solution $w$ to \eqref{yamabe0} and by the maximum principle, we know that $w>0$ in $M$ and $w(x)\to 1$ at infinity. So the metric $\breve{g}=w^{\frac{4}{n-2}}g_0$ is a complete Yamabe metric with zero scalar curvature, which is also quasi-isometric equivalent to the given metric $g_0$. We show that if the initial metric $g(0)$ is near to the metric $\breve{g}$, we have the global Yamabe flow which converges to the metric $\breve{g}$.
For the precise statement of this result, we refer to Theorem \ref{yamabe} in section \ref{sect2}. The Yamabe flow on $R^n$ belongs to the large class of porous-media equations, one may refer to \cite{Ar} and \cite{DK} for references and related results. The concept of a fine solution to local Yamabe flow has first been proposed in the work \cite{CZ}. The Yamabe flow on a hyperbolic space has a different interesting character as posed in the work \cite{S}. One may consult the paper \cite{MW} for the numerical method for the Yamabe metrics and Yamabe flows in dimension two.

We shall show that on an n-dimensional complete Riemannian manifold $(M,g_0)$ with non-negative Ricci curvature, $n\geq 3$, the Yamabe flow enjoys the local $L^1$-stability property. The precise result is stated in section \ref{sect3}.

We now consider the Yamabe problem on a  complete Riemannian manifold  with non-negative scalar curvature $R_0$.
Almost under the same condition as the Yamabe flow above, we can prove the following result about Yamabe metric with zero scalar curvature.
\begin{Thm}\label{zero}
Assume that $(M,g_0)$ is a complete Riemannian manifold  with non-negative scalar curvature $R_0$ such that the Poisson equation
$$
-\frac{4(n-1)}{n-2} \Delta_{g_0} u=R_0, \ \ in \ \ M, \ \ u(x)\to 0\ \ as \ d_{g_0}(x,p)\to \infty,
$$
has a positive solution $u$. Then there is a complete Yamabe metric $g$, conformal to the metric $g_0$, with zero scalar curvature.
\end{Thm}
We shall prove a very general result which will imply Theorem \ref{zero} above. See section \ref{sect1}.
Note that our result above is closely related to question (c) in page 234 in the work of J.Kazdan \cite{K} and the question says that if $M$ has complete metrics $g_+$ and $g_-$ with positive (respectively negative) scalar curvature, is there one with zero scalar curvature?
One may refer to the works \cite{A}, \cite{BE}, \cite{KW}, and \cite{SY}  for general discussions about scalar curvature problems.
Our main results are Theorems\ref{zero1}, Theorem \ref{zero2}, Theorem \ref{yamabe} and Theorem \ref{unique} in the latter sections.

The plan of this paper is below. In section \ref{sect1}, after a definition of property M on the complete Riemannian manifold with non-negative scalar curvature, we construct a a complete Yamabe metric on $(M^n,g_0)$. We study the fine solution to the Yamabe flow on the complete Riemannian manifold with non-negative Ricci curvature in section \ref{sect2}. In section \ref{sect3}, we consider $L^1$ stability of the Yamabe flow on the complete Riemannian manifold with non-negative scalar curvature. In the last section, we give some comments about the radially symmetric Yamabe metrics on model spaces.

\section{complete Yamabe metric with zero scalar curvature}\label{sect1}

To make our idea more transparent, we start from a concrete case.
We first give a definition. Let $(M,g_0)$ be a complete Riemannian manifold  with its scalar curvature $R_0=R(g_0)$.

\begin{Def}\label{def1}: We say the manifold $(M,g_0)$ satisfies \emph{property (M)}
if there is a smooth positive solution $v>0$ to the Poisson equation
\begin{equation}\label{poisson}
-\frac{4(n-1)}{n-2} \Delta_{g_0} u+R_0 u=R_0, \ \ in \ \ M,
\end{equation}
and
$$
u(x)\to 0, \ \ as\ \ x\to \infty.
$$
\end{Def}

Assume that $(M,g_0)$ is a complete Riemannian manifold  with non-negative scalar curvature $R_0$ and with \emph{property (M)}.
We want to show that on such a manifold, there is a complete Yamabe metric with zero scalar curvature.
Let $v$ be the solution to \eqref{poisson}.
We may let $d=\sup_M v(x)>0$. Note that $u=1$ is a solution to \eqref{poisson}.
Let $w(x)=1-v(x)$. Then
$$
-\frac{4(n-1)}{n-2} \Delta_{g_0} w+R_0 w=0.
$$
Since $w(x)\to 1$ as $x\to\infty$, by the maximum principle we know that $w(x)> 0$ on $M$ (
this implies that $d\leq 1$ on $M$). Note that the metric $w^{\frac{4}{n-2}}g_0 (\leq g_0)$ is the Yamabe metric with zero scalar curvature.

We now claim that for $(M,g_0)$ being a complete Riemannian manifold with non-negative scalar curvature $R_0$ and with the property that the Poisson equation
$$
-\frac{4(n-1)}{n-2} \Delta_{g_0} u=R_0, \ \ in \ \ M, \ \ u(x)\to 0\ \ as \ d(x,p)\to \infty
$$
has a positive solution $u$, there is a positive solution $v$ to \eqref{poisson} with decay at infinity. That is to say, property (M) is true on such a manifold.

In fact, since $u$ satisfies
$$
-\frac{4(n-1)}{n-2} \Delta_{g_0} u+R_0u-R_0=R_0 u\geq 0, \ \ in \ \ M,
$$
we know that $u$ is a super solution to \eqref{poisson}. Let $B_R=B_R(p)$ be the ball in $(M,g_0)$. We now solve, for $R>0$ large
\begin{equation}\label{ball}
-\frac{4(n-1)}{n-2} \Delta_{g_0} u+R_0u=R_0, \ u>0, \ in \ \ B_R,
\end{equation}
with $u=0$ on the boundary $\partial B_R$. We may get the positive solution by the use of Green function of the conformal Laplacian operator on $B_R$ with zero boundary condition. Let $u_R$ be the solution to \eqref{ball} on the ball $B_R$. Clearly, by the maximum principle, $u_R\leq u_{R'}\leq u$ on $B_R$ for $R'>R$. Hence, using the standard elliptic theory, we may get the $C^{3}_{loc}$-limit
$$
\bar{u}(x)=\lim_{R\to \infty} u_R(x), \ \ x\in M,
$$
which is a positive solution to \eqref{poisson} with the required property as in Property (M). Such a solution is unique via a use of maximum principle \cite{A}.
This proves the Claim.

We then conclude the below.

\begin{Thm}\label{zero1}
Assume that $(M,g_0)$ is a complete Riemannian manifold  with non-negative scalar curvature $R_0$ such that the property (M) is true on $M$. Then there is a complete Yamabe metric $g$, conformal to the metric $g_0$, with zero scalar curvature.
\end{Thm}

As the easy consequence of this result and  the argument above we can see Theorem \ref{zero} is true.

We remark that \emph{property (M)} is always true on asymptotic flat manifold $(M,g_0)$ with non-trivial scalar curvature $R_0\geq 0$ and
$R_0(x)=O(d(x,p)^{-\tau})$ near infinity for some $\tau>2$. This can be verified by the use of the Green's function.

\begin{Def}\label{def2}: We say the manifold $(M,g_0)$ satisfies \emph{property (H)}
if the following two conditions hold true on $M$.

(1). There is a smooth positive function $\phi>0$ such that the metric $\phi^{\frac{4}{n-2}}g_0$ is a complete metric with non-negative scalar curvature, i.e.,
$$
f:=L\phi=-\frac{4(n-1)}{n-2} \Delta_{g_0} \phi+R_0 \phi\geq 0, \ \ in \ \ M.
$$

(2). There is a positive solution $\bar{u}$ to the Poisson equation
\begin{equation}\label{poisson2}
-\frac{4(n-1)}{n-2} \Delta_{g_0} \bar{u}=f, \ \ in \ \ M
\end{equation}
such that for some constant $\theta\in (0,1)$,
$$
\frac{\bar{u}(x)}{\phi(x)} \leq \theta\in (0,1), \ \ as\ \ x\to \infty.
$$
\end{Def}

Assume that $(M,g_0)$ is a complete Riemannian manifold  with non-negative scalar curvature $R_0$ and with \emph{property (H)}.
We now show that on $(M,g_0)$, there is a complete Yamabe metric with zero scalar curvature.
Let $\bar{u}$ be the solution to \eqref{poisson2}. Note that
$$
L\bar{u}-f=R_0 \bar{u}\geq 0, \ \ in \ M.
$$
That is to say, $\bar{u}$ is a super-solution the equation
\begin{equation}\label{ball0}
Lv=f, \ \ in \ M.
\end{equation}

For $R>0$ large, we solve
\begin{equation}\label{ball2}
-\frac{4(n-1)}{n-2} \Delta_{g_0} u+R_0u=f, \ u>0, \ in \ \ B_R,
\end{equation}
with $u=0$ on the boundary $\partial B_R$. We then get the positive solution $u_R$ by the use of Green function of the conformal Laplacian operator on $B_R$ with zero boundary condition. By the maximum principle, $u_R\leq u_{R'}\leq \bar{u}$ on $B_R$ for $R'>R$. Hence, using the standard elliptic theory, we may get the $C^{3}_{loc}$-limit
$$
{u}_\infty(x)=\lim_{R\to \infty} u_R(x)\leq \bar{u}(x), \ \ x\in M
$$
which solves the equation \eqref{ball0}. Since $f=L\phi$, we then have, for $w(x):=\phi(x)-{u}_\infty(x)$ on $M$,
$Lw=0$ on $M$ and
$w(x)\geq (1-\theta)\phi(x)$ near infinity. By the maximum principle \cite{A}, we have $w(x)>0$ on $M$ and then the metric
$w^{\frac{4}{n-2}}g_0$ is a complete Yamabe metric with zero scalar curvature.

This proves the following result.
\begin{Thm}\label{zero2}
Assume that $(M,g_0)$ is a complete Riemannian manifold  with non-negative scalar curvature $R_0$ such that the property (H) is true on $M$. Then there is a complete Yamabe metric $g$, conformal to the metric $g_0$, with zero scalar curvature.
\end{Thm}

One may find more examples about Yamabe metrics with zero scalar curvature on the radial symmetric space models in section \ref{sect4}.

\section{global existence of fine solutions to the Yamabe flow}\label{sect2}
In this section, we assume that $(M,g_0)$ is a complete Riemannian manifold with Ricci curvature bounded from below and with non-negative scalar curvature. Assume also that the Poisson equation \eqref{poisson} has a positive solution $v(x)$, which is decay to zero at infinity.

Recall that the Yamabe flow is the metric family $(g(t))$ satisfying the evolution equation
$$
\partial_t g= -R(g)g, \ \  g=g(t)=u(t)^{\frac{4}{n-2}}g_0,
$$
with the given initial metric $g(0)$.
We denote by $u=u(t)$ and $g(t)=u^{\frac{4}{n-2}} g_0$.
We remark that we only assume the initial metric $g(0)$ is quasi isometric equivalent to the metric $g_0$.  More precisely, for the initial metric $g(0)$, we
assume that the metric $g(0)=u_0^{\frac{4}{n-2}}g_0$ is in the class
$$\mathbf{A}_C=\{u\in C^2(M); \ \ |u(x)-1|\leq Cv(x) \ in \ M\}$$
 where $C\geq 1$ is some uniform constant.

Note that
$$
R(g)=u^{-\frac{n+2}{n-2}}(-\frac{4(n-1)}{n-2} \Delta_{g_0} u+R_0 u), \ \  R_0=R(g_0).
$$
Then up to a scale change, we may write the Yamabe evolution equation as
\begin{equation}\label{yamabeflow1}
\partial_t (u^{\frac{n+2}{n-2}})=\frac{4(n-1)}{n-2} \Delta_{g_0} u-R_0 u,
\end{equation}
which can also be written as
\begin{equation}\label{yamabeflow2}
\frac{n+2}{n-2}u^{\frac{4}{n-2}}\partial_t u=\frac{4(n-1)}{n-2} \Delta_{g_0} u-R_0 u.
\end{equation}

The existence of the Yamabe flow with the given initial data $u_0:=u(0)$ can be constructed by the domain exhaustion method \cite{M2}.
Note that
$$
\frac{n+2}{n-2}u^{\frac{4}{n-2}}\partial_t (u-1)=\frac{4(n-1)}{n-2} \Delta_{g_0} (u-1)-R_0 (u-1)-R_0.
$$
Recall that for $C\geq 1$,
$$
-\frac{4(n-1)}{n-2} \Delta_{g_0} Cv+R_0 Cv=CR_0.
$$
Then we have
$$
-\frac{4(n-1)}{n-2} \Delta_{g_0}(-1+ Cv)+R_0 (-1+Cv)=(C-1)R_0, \ \ in \ M,
$$
which implies that $1-Cv$ is a subsolution to the Yamabe equation $Lu=0$ in $M$.
By this we have from \eqref{yamabeflow2},
$$
\frac{n+2}{n-2}u^{\frac{4}{n-2}}\partial_t (u-1+Cv)=\frac{4(n-1)}{n-2} \Delta_{g_0} (u-1+Cv)-R_0 (u-1+Cv)+(C-1)R_0.
$$
Then we have on $M$,
$$
\frac{n+2}{n-2}u^{\frac{4}{n-2}}\partial_t (u-1+Cv)\geq \frac{4(n-1)}{n-2} \Delta_{g_0} (u-1+Cv)-R_0 (u-1+Cv).
$$
By the maximum principle we know that
$$
u-1+Cv\geq 0, \ \ on  \ \ M
$$
provided the initial data satisfies $u_0-1+Cv\geq 0$ on $M$.

Similarly, we have, for $C\geq 1$,
$$
u-1-Cv\leq 0, \ \ on  \ \ M
$$
provided the initial data satisfies $u_0-1-Cv\leq 0$ on $M$. As for the convergent part  of the global Yamabe flow, we may just invoke the argument of Theorem 6 in \cite{M}.

We then conclude below.
\begin{Thm}\label{yamabe}
Assume that $(M,g_0)$ is a complete Riemannian manifold with Ricci curvature bounded from below and with non-negative scalar curvature. Assume that the Poisson equation \eqref{poisson} has a positive solution $v(x)$, which is decay to zero at infinity. Assume that the metric $g(0)=u_0^{\frac{4}{n-2}}g_0$ is in the class
$$\mathbf{A}_C=\{u\in C^2(M); \ \ |u(x)-1|\leq Cv(x) \ in \ M\}$$
 where $C\geq 1$ is some uniform constant. Then there exists a global Yamabe flow \eqref{yamabeflow1}  in the class $\mathbf{A}_C$ with the initial metric $g(0)$ and the flow converges to the Yamabe metric with zero scalar curvature.
\end{Thm}

\section{$L^1$ stability of the Yamabe flow}\label{sect3}

In this section, we study the local $L^1$ stability of the Yamabe flow \eqref{yamabeflow1} on the complete Riemannian manifold $(M,g_0)$ with non-negative Ricci curvature.
We prefer to rewrite the Yamabe flow equation as in the porous-media equation type
\begin{equation}\label{yamabeflow3}
\partial_t u=\frac{4(n-1)}{n-2} \Delta_{g_0} u^m-R_0 u^m,
\end{equation}
with $m=\frac{n-2}{n+2}$. Here, of course, the meaning of $u$ is different from the previous sections. Let $\alpha=\frac{n+2}{4}$.
Assume that we are given two solutions $u$ and $v$ to \eqref{yamabeflow3}.
By the Kato inequality we know that
$$
\Delta_{g_0} |u^m-v^m|\geq sign(u-v) \Delta_{g_0} (u^m-v^m).
$$
Then from the relation
$$
\partial_t (u-v)=\frac{4(n-1)}{n-2} \Delta_{g_0} (u^m-v^m)-R_0 (u^m-v^m).
$$
we get that
\begin{equation}\label{subharmonic}
\partial_t |u-v|\leq \frac{4(n-1)}{n-2} \Delta_{g_0} |u^m-v^m|-R_0 |u^m-v^m)|
\end{equation}
in the distributional sense.

We now consider the general case where the mean value property for any non-negative sub-harmonic function is true.
For any $\psi \in C_0^3(M)$, $0\leq \psi\leq 1$, we have
$$
\partial_t \int |u-v|\psi \leq \frac{4(n-1)}{n-2}  \int |\Delta_{g_0} \psi||u^m-v^m|.
$$
By Holder inequality we have
$$
\int |\Delta_{g_0} \psi||u^m-v^m|\leq \int |\Delta_{g_0} \psi||u-v|^m\leq C(\psi)[\int |\psi||u-v|]^m,
$$
where
$$
C(\psi)=\int |\Delta_{g_0} \psi|^{\frac{n+2}{4}}|\psi|^{-\frac{n+2}{4}m}]^{1-m}.
$$
Then we have
$$
\partial_t \int |u-v|\psi \leq\frac{4(n-1)}{n-2}C(\psi)[\int |\psi||u-v|]^m.
$$
Integrating the relation above, we get for any $t,s\geq 0$,
\begin{equation}\label{stable}
[ \int |u-v|\psi(t)]^{1-m} \leq  [\int |u-v|\psi (s)]^{1-m}+ (1-m)\frac{4(n-1)}{n-2} C(\psi) |t-s|.
\end{equation}
We call the relation \eqref{stable} the local $L^1$ stability inequality.
Choose $\psi=\psi_0^k$, $k>2\alpha$ large, $\psi_0$ is the cut-off function on $B_{R}(p)$ for $R>0$. Then
\begin{equation}\label{constant}
C(\psi)\leq C(m,n)[\frac{V(R)}{R^{2\alpha}}]^{1-m}, \ \
\end{equation}
where $V(R)=vol(B_R(p))$.

Assume now that $u(0)=v(0)$. Then we have
\begin{equation}\label{mean}
\int_{B_R(p)} |u-v| \leq C t^{\frac{1}{1-m}}\frac{V(R)}{R^{2\alpha}}.
\end{equation}
We let
$$
w(t,x)=\int_0^t |u^m-v^m|(s,x)ds.
$$
Since, for $a>b>0$,
$$
|a^m-b^m|\leq 2^{1-m}|a-b|^m,
$$
we have
$$
w(t,x)\leq 2^{1-m} \int_0^t|u-v|^m(s,x)ds.
$$
By \eqref{subharmonic} we know that
$$
0\leq |u(t,x)-v(t,x)|\leq \frac{4(n-1)}{n-2} \Delta_{g_0} w(t,x).
$$
By the mean value property (see \cite{LS} or Schoen-Yau's book \cite{SY}) we know that for any $R>0$,
\begin{equation}\label{mean2}
w(t,p)\leq \frac{C}{V(R)}\int_{B_R(p)} w(t,x).
\end{equation}
Note that
$$
\int_{B_R(p)} w(t,x)\leq C\int_0^t \int_{B_R(p)} |u-v|^m\leq CV(R)^{1-m}\int_0^t ds(\int_{B_R(p)} |u-v|)^m.
$$
That is to say,
$$
\frac{1}{V(R)}\int_{B_R(p)} w(t,x)\leq C\int_0^t ds(\frac{1}{V(R)}\int_{B_R(p)} |u-v|)^m.
$$
By \eqref{mean} we get
$$
\int_0^t ds(\frac{1}{V(R)}\int_{B_R(p)} |u-v|)^m\leq \int_0^t ds s^{\frac{m}{1-m}}(\frac{1}{R^{2\alpha}})^m=C t^{\frac{1}{1-m}}(\frac{1}{R^{2\alpha}})^m
$$
and
$$
\frac{1}{V(R)}\int_{B_R(p)} w(t,x)\leq C t^{\frac{1}{1-m}}(\frac{1}{R^{2\alpha}})^m.
$$

Putting this back to \eqref{mean2} we know that for some $C(t)>0$,
$$
w(t,p)\leq \frac{C(t)}{R^{2m\alpha}}\to 0, \ \ as\ \ R\to\infty.
$$
This implies the uniqueness desired.

\begin{Thm}\label{unique} Assume that $(M,g_0)$ is a complete Riemannian manifold with non-negative Ricci curvature.
The domain exhaustion solutions to the Yamabe flow satisfy the local $L^1$ stability inequality \eqref{stable}  and the exhaustion solution to the Yamabe flow with given initial data is
 unique.
\end{Thm}

In the proof above, we have used an idea from the paper \cite{HP}. We may make some estimate above more precise in the case when $(M,g_0)$ is a complete Riemannian manifold with non-negative Ricci curvature.
By Bishop volume comparison theorem \cite{SY} we have
$$
V(R)=vol(B_R(p))\leq C_n R^n
$$ on the manifold $(M,g_0)$ with non-negative Ricci curvature. In this case, we have the estimate of the constant in \eqref{constant} that
$$
C(\psi)\leq C(m,n) R^{(-\frac{n}{2}+1)(1-m)}.
$$
Then we have the local stability bound that
$$
[ \int |u-v|\psi(t)]^{1-m} \leq  [\int |u-v|\psi (0)]^{1-m}+C(m,n)tR^{(-\frac{n}{2}+1)(1-m)}.
$$
This local uniform bound may be useful in the local in time stability of the Yamabe flow.

\section{comments about the Yamabe metric on radially symmetric spaces}\label{sect4}

In this section we briefly discuss the radially symmetric Yamabe metrics on a model space, which is a radially symmetric space $M,g_0)$ \cite{BK},
 where
 $$
 M=N\times R_+, \ \ g_0=dr^2+ f^2(r)d\sigma^2, \ \ f(0)=0, \ f'(0)=1,
 $$
 and $(N, d\sigma^2)$ is a compact Riemannian manifold with positive constant curvature. So we may assume that $M=R^n$ and $N=S^{n-1}$ the sphere. We look for radial Yamabe metric and we shall reduce the problem into an ODE, which is of Riccati type ODE and may be approachable by classical ODE methods. Our discussions about the related ODE are abstract in principle.

 Note that
 $$
 f''(r)=-K(r) f(r),
 $$
where $K(r)=K(|x|)$ is the radial curvature of $(R^n,g_0)$. We let
$$
K_1(r)=\frac{1-f_r^2}{f^2}
$$
which is the tangent curvature to the spheres.
Note that the Laplacian operator is
$$
\partial_r^2+\frac{(n-1)f'(r)}{f(r)}\partial_r+\frac{1}{f^2(r)}\Delta_\sigma.
$$
The Ricci curvature is
$$
Rc=(n-1)K dr^2+(K+(n-2)K_1)d\sigma^2
$$
and the scalar curvature of $g_0$ is
$$
R_0=2(n-1)K+(n-1)(n-2)K_1
$$
respectively.

We want to find a metric $g_v=v^{\frac{4}{n-2}}g_0$, where $v=v(r)$, such that
$$
-\frac{4(n-1)}{n-2} \Delta_{g_0} v+R_0 v=0, \ \ on \ \ M.
$$
That is to say, $v=v(r)$ satisfies the ODE
$$
-(v''+\frac{(n-1)f'(r)}{f(r)}v')+\frac{n-2}{4(n-1)}R_0v=0, \ \ r\geq 0
$$
with the initial datum
$$
v(0)=v_0>0, v'(0)=v_1.
$$
Define
$$
L_1=-f^{\frac{n-1}{2}}\Delta_{g_0} f^{-\frac{n-1}{2}}.
$$
Then
$$
L_1=-\frac{d^2}{dr^2}+Q(r),
$$
where
$$
Q(r)=\frac{(n-1)(n-3)}{4}\frac{f'(r)^2}{f(r)^2}-\frac{n-1}{2}K(r).
$$

Define $V=f^{-\frac{n-1}{2}}v$. Then we have
$$
L_1V+\frac{n-2}{4(n-1)}R_0V=0, \ \ r\geq 0.
$$
In explicit form, we have
$$
-V''+[Q(r)+\frac{n-2}{4(n-1)}R_0]V=0,
$$
which may be solved by the classical ODE method.

Assume that there are constant $c,d>0$ and a smooth function $a(r)$ such that either (1)
$$
a'=-a^2+[Q(r)+\frac{n-2}{4(n-1)}R_0], \ \ on  \ [0,\infty)
$$
and with the function
$$
f(r)^{\frac{n-1}{2}}exp (\int_0^r a(t)dt) \to c>0, \ \ as \ \ r\to \infty;
$$
or (2)
$$
a'=a^2+[Q(r)+\frac{n-2}{4(n-1)}R_0], \ \ on  \ [0,\infty)
$$
and with the function
$$
f(r)^{\frac{n-1}{2}}exp (-\int_0^r a(t)dt) \to d>0, \ \ as \ \ r\to \infty.
$$

In case (1), we have
$$
V''-[Q(r)+\frac{n-2}{4(n-1)}R_0]V=(\frac{d}{dr}+a)(\frac{d}{dr}-a)V=0.
$$
We may solve the first order ODE
$$
V'-aV=0,
$$
which gives the positive solution
$$
V(r)=V(0) exp(\int_0^r a(t)dt).
$$

In case (2), we have
$$
V''-[Q(r)+\frac{n-2}{4(n-1)}R_0]V=(\frac{d}{dr}-a)(\frac{d}{dr}+a)V=0.
$$
We may solve the first order ODE
$$
V'+aV=0,
$$
which gives the positive solution
$$
V(r)=V(0) exp(-\int_0^r a(t)dt).
$$

Recall that $v(r)=f(r)^{\frac{n-1}{2}} V(r)$.
In either case, we can get a complete Yamabe metric with zero scalar curvature.
We may formulate the above discussions into an abstract theorem, but we omit the detail.

\end{document}